\documentclass[a4paper,11pt]{article}
\usepackage{amsfonts}
\usepackage{amssymb}
\usepackage{amsmath}
\usepackage{graphicx}

\usepackage{amsthm}
\usepackage{cite}

\def\figurename{Figure}
\makeatletter
\renewcommand{\fnum@figure}[1]{\figurename~\thefigure}

\newtheorem{theorem}{Theorem}[section]
\newtheorem{property}{Property}[section]
\newtheorem{proposition}{Proposition}[section]

\numberwithin{equation}{section}

\begin{document}

\title{\textbf{Two Regular Polygons \\ with a Shared Vertex}}
\author{\textbf{Mamuka Meskhishvili}}

\date{}
\maketitle

\begin{abstract}
For two non-congruent regular polygons of the same type, the method of finding the points in the plane at the equal distances to the vertices, is established. The existence of two points with this property is proved for two polygons with a shared vertex. For one of them, it is proved that it satisfies the Bottema theorem conditions and based on this, the generalized Bottema theorem for any two regular polygons is given.

\vskip1em \noindent \textbf{2010 AMS Classification:}
    51M15, 51M20, 51M35.

\vskip1em \noindent \textbf{Keywords and phrases:}
    Bottema theorem, cyclic averages, Regular polygon, Common vertex.

\end{abstract}

\bigskip
\bigskip

\section{Introduction}
\label{sec:1}

\bigskip

For a regular $n$-sided polygon $A_1A_2\cdots A_n$ and an arbitrary point $M$ in the plane of the polygon, the distances from $M$ to the vertices $A_1,A_2,\dots,A_n$ satisfy \cite{1, 2}:
\begin{equation}\label{eq:1}
    \sum_{i=1}^n d_i^{2m}=n\Bigg[(R^2+L^2)^m+\sum_{k=1}^{\lfloor\frac{m}{2}\rfloor} \binom{m}{2k}\binom{2k}{k}R^{2k}L^{2k}(R^2+L^2)^{m-2k}\Bigg],
\end{equation}
where $m=1,\dots,n-1$; $R$ is the radius of the circumcircle $\Omega$ and $L$ is the distance between $M$ and the centroid $O$ of the regular polygon.

Let us the second $n$-sided polygon $B_1B_2\cdots B_n$ is given in the plane. The distance from the point $M$ to the vertices $B_1,B_2,\dots,B_n$ denote by $t_1,t_2,\dots,t_n$. Therefore, for given two $n$-sided regular polygons and the point $M$, we have two sets of distances:
$$  \{d_i\} \;\;\text{and}\;\; \{t_i\}.     $$
Are there points in the plane, which have the same set of these distances? This problem is investigated in the present paper.

Denote by $R_2$ and $L_2$ the radius of the circumcircle $\Omega_2$ and the distance between $M$ and centroid $O_2$ of the regular polygon $B_1B_2\cdots B_n$. Equalize the right sides of \eqref{eq:1}, we get:
\begin{multline*}
    (R_1^2+L_1^2)^m+\sum_{k=1}^{\lfloor\frac{m}{2}\rfloor} \binom{m}{2k}\binom{2k}{k}R_1^{2k}L_1^{2k}(R_1^2+L_1^2)^{m-2k} \\
    =(R_2^2+L_2^2)^m+\sum_{k=1}^{\lfloor\frac{m}{2}\rfloor} \binom{m}{2k}\binom{2k}{k}R_2^{2k}L_2^{2k}(R_2^2+L_2^2)^{m-2k}.
\end{multline*}
For any $n$-gons, the first two relations $(m=1,2)$ are:
\begin{align*}
    R_1^2+L_1^2 & =R_2^2+L_2^2, \\ 
    (R_1^2+L_1^2)^2+2R_1^2L_1^2 & =(R_2^2+L_2^2)^2+2R_2^2L_2^2.    
\end{align*}

So, we obtain two cases:
\begin{align*}
    \text{congruent case -- } &\; R_1=R_2 \;\;\text{and}\;\; L_1=L_2, \\
    \text{non-congruent case -- } &\; R_1=L_2 \;\;\text{and}\;\; L_1=R_2.
\end{align*}

Equalize the left sides of \eqref{eq:1}, we get the $n-1$ relations for distances:
\begin{equation}\label{eq:*}
\begin{aligned}
    d_1^2+d_2^2+\cdots+d_n^2 & =t_1^2+t_2^2+\cdots+t_n^2, \\
    d_1^4+d_2^4+\cdots+d_n^4 & =t_1^4+t_2^4+\cdots+t_n^4, \\
    & \vdots \\
    d_1^{2(n-1)}+d_2^{2(n-1)}+\cdots+d_n^{2(n-1)} & =t_1^{2(n-1)}+t_2^{2(n-1)}+\cdots+t_n^{2(n-1)}.
\end{aligned}       \tag{$*$}
\end{equation}

Let's consider the cases separately.

\bigskip
\bigskip

\section{Congruent regular polygons}
\label{sec:2}

\bigskip

Congruent $n$-gons case is divided into two subcases $O_1=O_2$ and $O_1\neq O_2$.

\begin{figure}[h]
\centerline{\includegraphics[width=7cm]       
    {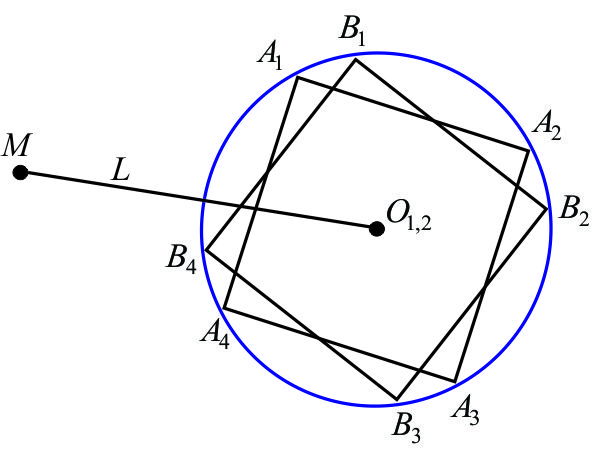} }
\caption{}
\label{fig:1}
\end{figure}

\begin{figure}[h]
\centerline{\includegraphics[width=9.3cm]       
    {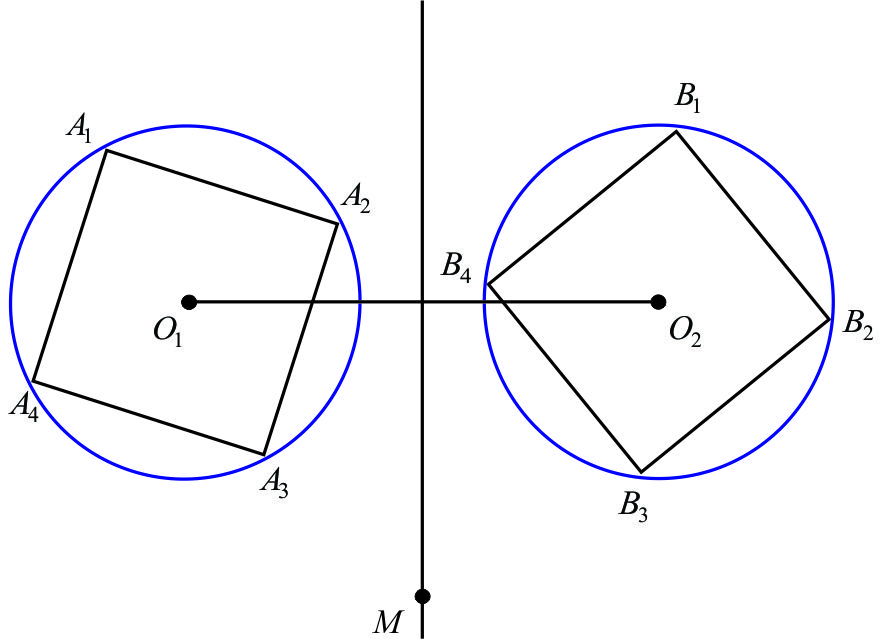} }
\caption{}
\label{fig:2}
\end{figure}

If two polygons centroids coincide, we get two $n$-gons inscribed in the same circle and they differ only by rotation around the common centroid, see Fig.~\ref{fig:1} (for figures of any $n$-gons we use the squares in the Fig.~\ref{fig:1}--\ref{fig:6}).

For $M$, we can take any point in the plane.

\begin{proposition}[Rotational invariant]\label{prop:1}
If two congruent polygons are inscribed in the same circle, the distances from any point in the plane to the vertices of the polygons satisfy the system~\eqref{eq:*}.
\end{proposition}

If $O_1\neq O_2$, it is clear $M$ lies on the perpendicular bisector of the line segment $O_1O_2$, see Fig.~\ref{fig:2}.

For $M$, we can take any point on the perpendicular bisector.

\begin{proposition}[Reflection invariant]\label{prop:2}
If two congruent polygons are inscribed in symmetrical circles, the distances from any point on the axis of symmetry to the vertices of polygons satisfy the system \eqref{eq:*}.
\end{proposition}

\bigskip
\bigskip

\section{Non-congruent regular polygons}
\label{sec:3}

\bigskip

\begin{figure}[h]
\centerline{\includegraphics[width=7cm]       
    {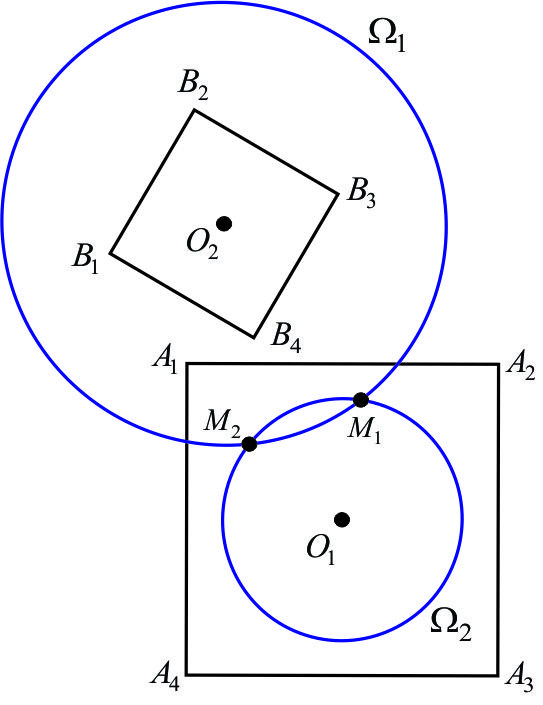} }
\caption{}
\label{fig:3}
\end{figure}

From the conditions of non-congruent case --
$$  R_1=L_2 \;\;\text{and}\;\; L_1=R_2      $$
follow -- the point $M$ is the intersection point of two circles:
$$  \Omega_1(O_2,R_1) \;\;\text{and}\;\; \Omega_2(O_1,R_2).     $$
It is clear, such point exists, if
\begin{equation}\label{eq:5}
    |R_1-R_2|\leq O_1O_2\leq R_1+R_2,
\end{equation}
and if the circles intersect each other there are two points $M_1$ and $M_2$ of such property see Fig.~\ref{fig:3}

For $M$, we can take only points $M_1$ and $M_2$.

\begin{proposition}[Necessary condition]\label{prop:3}
For two non-congruent $n$-gons with centroids -- $O_1$, $O_2$, circumradiii -- $R_1$, $R_2$ and circumcircles $\Omega_1$, $\Omega_2$, the only points of the intersection
$$  \{M_1,M_2\}=\Omega_1(O_2,R_1)\cap\Omega_2(O_1,R_2)      $$
satisfy the system \eqref{eq:*}.
\end{proposition}

Until now we consider all cases for which the distances from the point to the vertices of two regular $n$-gons satisfy the system \eqref{eq:*}, but for equality of the distances the system \eqref{eq:*} is only necessary condition. Let us establish the sufficient condition.

\bigskip
\bigskip

\section{Equalization of distances}
\label{sec:4}

\bigskip

Let us consider the distances from the given point $M_1$ to the vertices of the second regular polygon $B_1B_2\cdots B_n$ as variables:
$$  t_1,t_2,\dots,t_n.      $$

In the system \eqref{eq:*}, we have $n-1$ equations and $n$ variables in order to determine the variables uniquely we must eliminate one of them. It is possible by using Proposition \ref{prop:1}. The rotation does not change the system \eqref{eq:*} thus we can rotate the second polygon, so that one pair of distances will be equal to each other.

Explain this procedure by using Fig.~\ref{fig:3}. Construct circumscribe circles of the second polygon -- $\Omega_2$, whose center is $O_2$. From the point $M_1$ as the center draw the auxiliary circle with radius $M_1A_1$.

\begin{figure}[h]
\centerline{\includegraphics[width=7cm]       
    {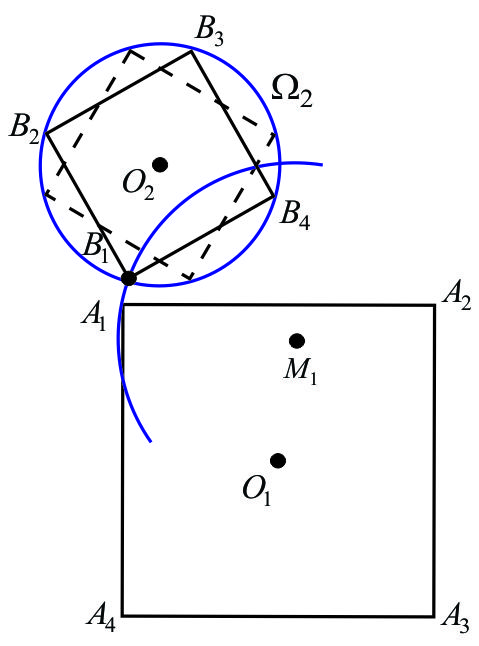} }
\caption{}
\label{fig:4}
\end{figure}

The intersection point of the auxiliary and $\Omega_2$ circles is the new position of the vertex $B_1$, see Fig.~\ref{fig:4}. For the new position of the polygon $B_1B_2\cdots B_n$, we obtain:
\begin{equation}\label{eq:**}
    t_1=d_1.        \tag{$**$}
\end{equation}

The condition \eqref{eq:**} and the system \eqref{eq:*}, give us a new system:
\begin{equation}\label{eq:***}
\begin{aligned}
    d_2^2+\cdots+d_n^2 & =t_2^2+\cdots+t_n^2, \\
    d_2^4+\cdots+d_n^4 & =t_2^4+\cdots+t_n^4, \\
    & \vdots \\
    d_2^{2(n-1)}+\cdots+d_n^{2(n-1)} & =t_2^{2(n-1)}+\cdots+t_n^{2(n-1)}.
\end{aligned}       \tag{$***$}
\end{equation}
Now the number of the variables equals to the number of the equations. By using Newton's Identities, elementary symmetric polynomials can be expressed in terms of power sums \cite{4, 6}. Because of \eqref{eq:***}, the power sums of $(d_2^2,\dots,d_n^2)$ equal to the power sums of $(t_2^2,\dots,t_n^2)$, so the corresponding elementary polynomials are the same, so
$$  d_2^2,\dots,d_n^2 \;\;\text{and}\;\; t_2^2,\dots,t_n^2      $$
are the roots of the same equation of degree $n-1$. Consequently $t_2,\dots,t_n$ are the permutation of the $d_2,\dots,d_n$, therefore:

\begin{proposition}[Sufficient condition]\label{prop:4.1}
If in the system \eqref{eq:***}, one pair of the distances is the same $(t_1=d_1)$, from that it follows the equality of the sets
$$  \{t_2,\dots,t_n\}=\{d_2,\dots,d_n\}.        $$
\end{proposition}

For congruent case, the rotation gives the coincidence (if $O_1=O_2$, see Fig.~\ref{fig:1}) and reflexion symmetry (if $O_1\neq O_2$, see Fig.~\ref{fig:2}).

For non-congruent polygons, we obtain:

\begin{theorem}\label{th:4.1}
If two non-congruent regular polygons are given in the plane, it is possible to fix both of them, so that there is the point $M$ (in the plane) at equal distances from the vertices of the polygons. The point $M$ satisfies:
\begin{enumerate}
\item[{\rm I.}] $M=\Omega_1(O_2,R_1)\cap \Omega_2(O_1,R_2)$; where $\Omega_1$, $\Omega_2$ -- are circumscribed circles, $R_1$, $R_2$ -- circumradii and $O_1$, $O_2$ -- centroids of the polygons.

\item[{\rm II.}] One pair of distances must be equal, which is possible by rotation of one polygon over its centroid.
\end{enumerate}
\end{theorem}

By convenient enumeration of the vertices, we get
$$  MA_k=MB_k, \;\;\text{where}\;\; k=1,\dots,n.        $$

\bigskip
\bigskip

\section{A shared vertex}
\label{sec:5}

\bigskip

If two $n$-gons have a shared vertex, the equality \eqref{eq:**} automatically holds for $M_1$ and $M_2$, so we do not need the rotation. For non-congruent case, we obtain -- two regular polygons and two points theorem.

\begin{theorem}\label{th:4}
If in the plane two regular non-congruent polygons of the same type have a shared vertex, there are two points in the plane separately having equal distances to the vertices of the polygons.
\end{theorem}

The point $M_1$ and $M_2$ satisfy:
$$  \{M_1,M_2\}=\Omega_1(O_2,R_1)\cap\Omega_2(O_1,R_2).     $$
In the Fig.~\ref{fig:5} is given the case of two squares with shared vertex -- $A_1$:

\begin{figure}[h]
\centerline{\includegraphics[width=10cm]       
    {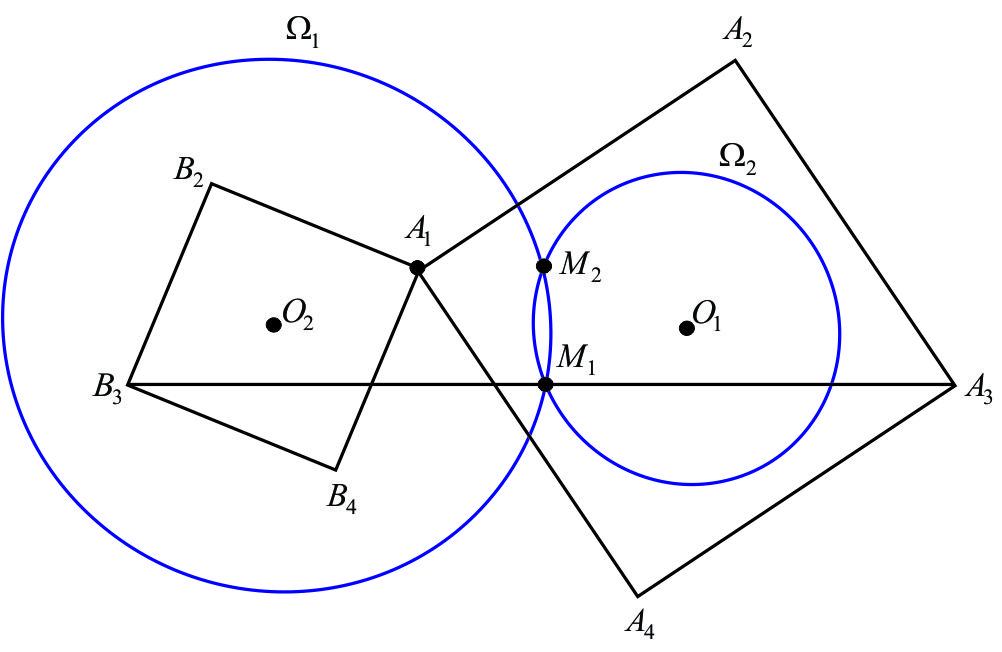} }
\caption{}
\label{fig:5}
\end{figure}

The same distances are:
\begin{equation}\label{eq:6}
    M_1A_2=M_1B_2, \quad M_1A_3=M_1B_3, \quad M_1A_4=M_1B_4
\end{equation}
and
\begin{equation}\label{eq:7}
    M_2A_2=M_2B_4, \quad M_2A_3=M_2B_3, \quad M_2A_4=M_2B_2.
\end{equation}

\begin{figure}[h]
\centerline{\includegraphics[width=10cm]       
    {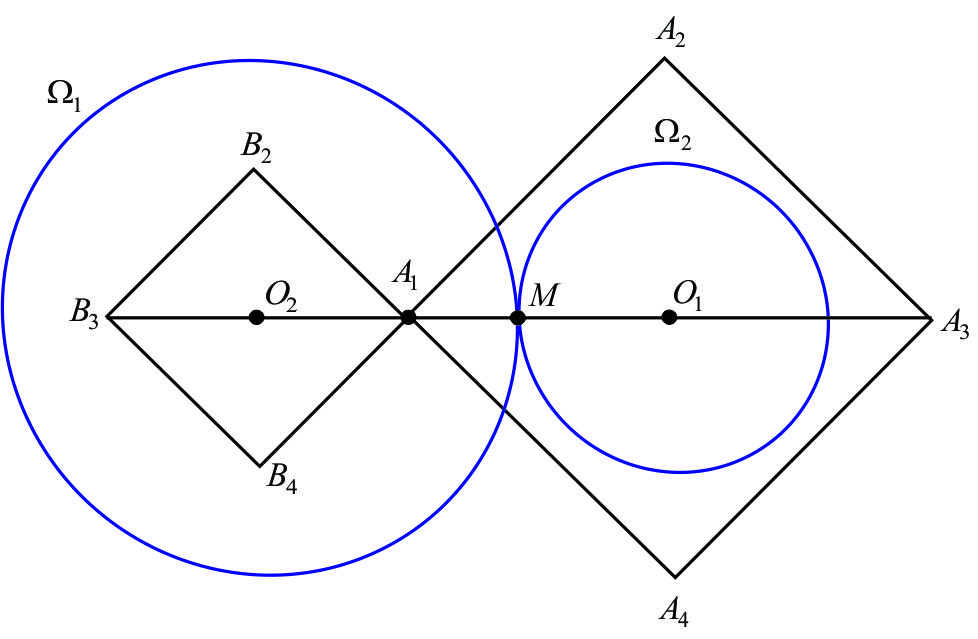} }
\caption{}
\label{fig:6}
\end{figure}

In case of the shared vertex, from the triangle -- $O_1A_1O_2$ the condition \eqref{eq:5} is always held, so generally there are two points. If the shared vertex and the centroids are collinear, there is only one point having equal distances to the vertices (see Fig.~\ref{fig:6})

\bigskip
\bigskip

\section{Corresponding equal distances}
\label{sec:6}

\bigskip

\begin{figure}[h]
\centerline{\includegraphics[width=10cm]       
    {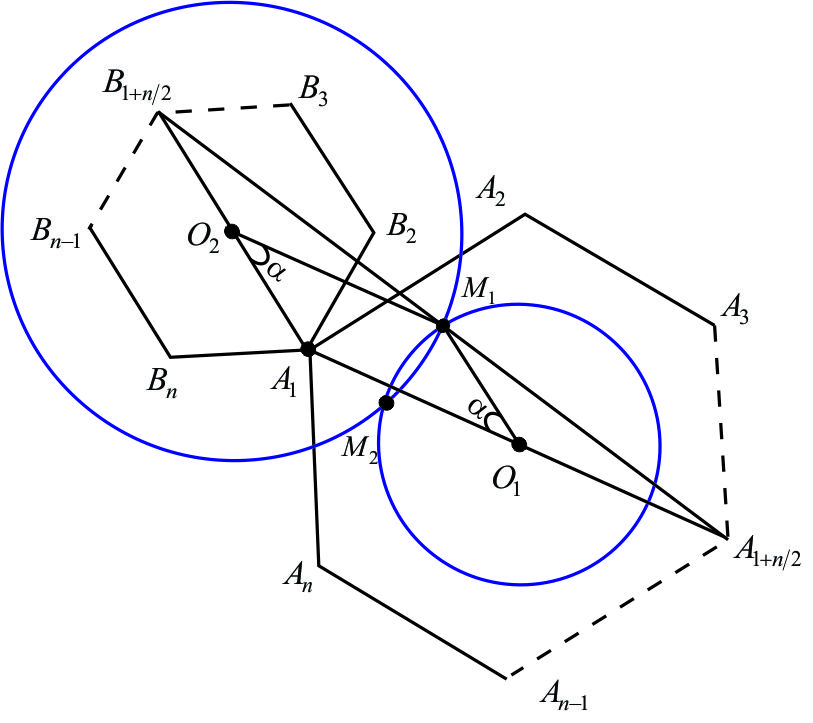} }
\caption{}
\label{fig:7}
\end{figure}

\begin{figure}[h]
\centerline{\includegraphics[width=10cm]       
    {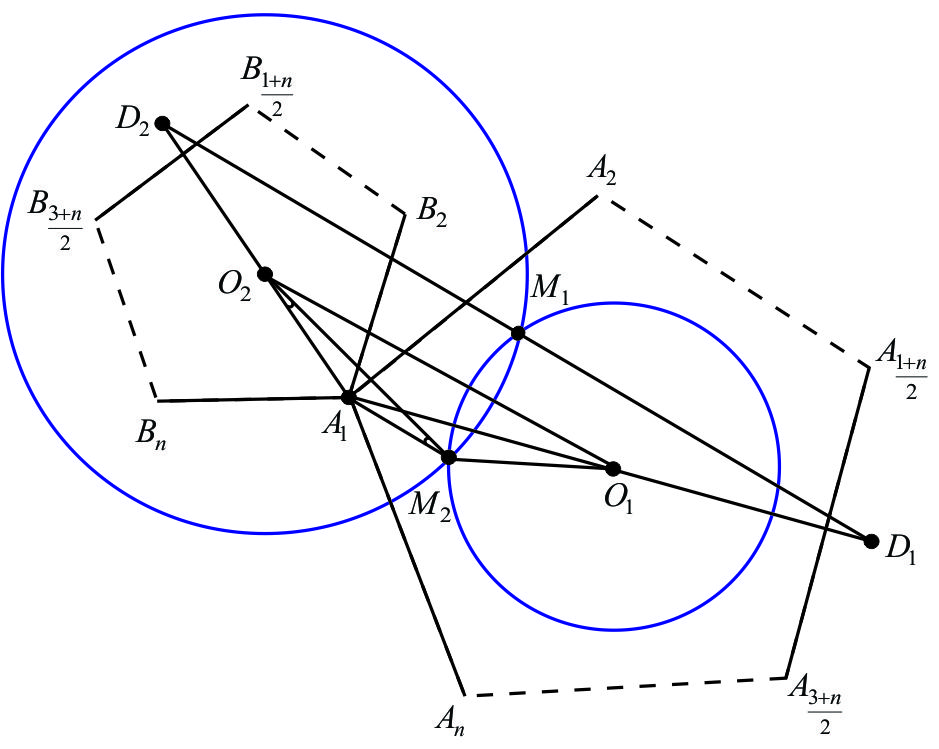} }
\caption{}
\label{fig:8}
\end{figure}

The equal distances \eqref{eq:6} and \eqref{eq:7}, for points $M_1$ and $M_2$ are different in the order. If we take the vertices of the first polygon in clockwise direction:
$$  A_2,\;A_3,\;A_4,;       $$
the vertices of the corresponding equal distances for the second polygon are:
\begin{enumerate}
\item[-] for $M_2$, in clockwise direction -- $B_4,\;B_3,\;B_2$;

\item[-] for $M_1$, in anticlockwise direction  -- $B_2,\;B_3,\;B_4$;
\end{enumerate}
i.e. the vertices of corresponding equal distances for points $M_1$ and $M_2$ are in opposite order . Is it true for any regular polygons of the same type? Let us consider two regular $n$-gons with the shared vertex $A_1$ (see Fig.~\ref{fig:7} and Fig.~\ref{fig:8}).

For the point $M_1$ in Fig.~\ref{fig:7}, we have
$$  \vartriangle\!M_1A_1O_2=\vartriangle\!M_1O_1A_1, \;\;\text{so}\;\; \angle\,A_1O_2M_1=\angle\,M_1O_1A_1\equiv\alpha.        $$
Then,
\begin{align*}
    M_1A_2^2 & =R_1^2+R_2^2-2R_1R_2\cos\Big(\frac{2\pi}{n}-\alpha\Big), \\
    M_1B_2^2 & =R_1^2+R_2^2-2R_1R_2\cos\Big(\frac{2\pi}{n}-\alpha\Big), \\
    &\qquad\qquad\qquad \vdots \\
    M_1A_k^2 & =R_1^2+R_2^2-2R_1R_2\cos\Big((k-1)\,\frac{2\pi}{n}-\alpha\Big), \\
     M_1B_k^2 & =R_1^2+R_2^2-2R_1R_2\cos\Big((k-1)\,\frac{2\pi}{n}-\alpha\Big).
\end{align*}
Therefore,
$$  M_1A_k=M_1B_k, \;\;\text{where}\;\; k=2,\dots,n,        $$
i.e. the vertices of the corresponding equal distances are in anticlockwise direction.

In the same manner for the point $M_2$, we have (see Fig.~\ref{fig:8}):
$$  \vartriangle\!M_2O_1A_1=\vartriangle\!M_2O_2A_1, \quad \angle\,M_2O_1A_1=\angle\,M_2O_2A_1\equiv\beta.        $$
Then,
\begin{align*}
    M_2A_2^2 & =R_1^2+R_2^2-2R_1R_2\cos\Big(\frac{2\pi}{n}+\beta\Big), \\
    M_2B_n^2 & =R_1^2+R_2^2-2R_1R_2\cos\Big(\frac{2\pi}{n}+\beta\Big), \\
    &\qquad\qquad\qquad \vdots \\
    M_2A_k^2 & =R_1^2+R_2^2-2R_1R_2\cos\Big((k-1)\,\frac{2\pi}{n}+\beta\Big), \\
     M_2B_{n+2-k}^2 & =R_1^2+R_2^2-2R_1R_2\cos\Big((k-1)\,\frac{2\pi}{n}+\beta\Big).
\end{align*}
Therefore,
$$  M_2A_k=M_2B_{n+2-k}, \;\;\text{where}\;\; k=2,\dots,n,        $$
i.e. the vertices of the corresponding equal distances are in anticlockwise direction.

We obtain:

\begin{theorem}\label{th:5}
There are two points $M_1$ and $M_2$ in the plane having the equal distances to the vertices of two regular polygons $A_1A_2\cdots A_n$ and    \linebreak   $A_1B_2\cdots B_n$ with a shared vertex $A_1$:
$$  M_1A_k=M_1B_k \;\;\text{and}\;\; M_2A_k=M_2B_{n+2-k}, \;\;\text{where}\;\;k=2,\dots,n.        $$
\end{theorem}

\bigskip
\bigskip

\section{Properties of $M_1$ and $M_2$}
\label{sec:7}

\bigskip

If the number $n$ of the sides of the polygons is even, then diametrically opposed points of the shared vertex $A_1$ are the vertices
$$  A_{1+\frac{n}{2}} \;\;\text{and}\;\; B_{1+\frac{n}{2}},     $$
see Fig.~\ref{fig:7}.

The quadrilateral $O_2M_1O_1M_2$ is parallelogram with sides -- $R_1$, $R_2$, then
$$  \angle B_{1+\frac{n}{2}}M_1A_{1+\frac{n}{2}}=\angle B_{1+\frac{n}{2}}M_1O_2+\angle O_2M_1O_1+
            \angle A_{1+\frac{n}{2}}M_1O_1=\pi,      $$
because of
$$  \angle B_{1+\frac{n}{2}}M_1O_2=\angle M_1A_{1+\frac{n}{2}}O_1 \;\;\text{and}\;\;
            \angle M_1A_{1+\frac{n}{2}}O_1+\angle A_{1+\frac{n}{2}}M_1O_1=\angle M_1O_1A_1.  $$

\begin{property}\label{pr:1}
If $n$ is even, the point $M_1$ is the midpoint of the line segment $A_{1+\frac{n}{2}}B_{1+\frac{n}{2}}$.
\end{property}

If the number $n$ is odd, diametrically opposed points $D_1$, $D_2$ of the $A_1$ are the midpoints of the arcs
$$  A_{\frac{1+n}{2}}A_{\frac{3+n}{2}} \;\;\text{and}\;\; B_{\frac{1+n}{2}}B_{\frac{3+n}{2}}\,,  $$
of the circles $\Omega_1(O_1,R_1)$, $\Omega_2(O_2,R_2)$, respectively (see Fig.~\ref{fig:8}).

In odd case we can ``double'' the number of the vertices and then use the Property~\ref{pr:1}. For even case
$$  D_1=A_{1+\frac{n}{2}} \;\;\text{and}\;\; D_2=B_{1+\frac{n}{2}},         $$
so for both cases, we have:

\begin{property}\label{pr:2}
The point $M_1$ is the midpoint of the line segment $D_1D_2$, where $D_1$, $D_2$ are diametrically opposed points of the shared vertex of the polygons.
\end{property}

The triangle $M_2D_1D_2$ is isosceles (Theorem \ref{th:5}, $k=1+\frac{n}{2}$), so
$$  M_1M_2\perp D_1D_2.         $$

In the quadrelateral $A_1M_2O_1O_2$ (see Fig.~\ref{fig:8})
$$  M_2O_1=O_2A_1 \;\;\text{and}\;\; A_1O_1=M_2O_2,     $$
so $A_1M_2O_1O_2$ is the isosceles trapezoid and
$$  A_1M_2\,/\!/\,O_2O_1.         $$
The line segment $O_1O_2$ is the midsegment of the triangle $D_1A_1D_2$, so we obtain:

\begin{property}\label{pr:3}
The point $M_2$ lies on the perpendicular bisector of $D_1D_2$ and $A_1M_2\,/\!/\,D_1D_2$.
\end{property}

\begin{property}\label{pr:4}
The line segment $M_1M_2$ is equal to the distance from the shared vertex to the line $D_1D_2$.
\end{property}

\bigskip
\bigskip

\section{Generalized Bottema theorem}
\label{sec:8}

\bigskip

The Bottema theorem concerns two squares, which have a common vertex. The theorem can be stated as follows (see Fig.~\ref{fig:5}) \cite{3, 5}:

\textit{In any given triangle $A_4A_1B_4$, construct two squares on two sides $A_1A_4$ and $A_1B_4$. The midpoint of the line segment that connects the vertices of squares opposite the common vertex $A_1$, is independent of the location of~$A_1$.}

Let us change two squares by two regular polygons of the same number of sides in the Bottema theorem:

\begin{theorem}\label{th:6}
In any given triangle $A_nA_1B_n$, construct two regular $n$-gons $A_1A_2\cdots A_n$ and $A_1B_2\cdots B_n$ on two sides $A_1A_n$ and $A_1B_n$. Take the points $D_1$, $D_2$ on the circumcircles of the polygons, which are diametrically opposed of the common vertex $A_1$. Then, the midpoint of the line segment $D_1D_2$ is independent of the location of $A_1$.
\end{theorem}

\vskip+0.2cm
\noindent \textit{Proof.}
In our notations (see Fig.~\ref{fig:7} and Fig.~\ref{fig:8}) the midpoint is $M_1$, and from Theorem \ref{th:5}:
\begin{equation*}
    M_1A_k=M_1B_k, \;\;\text{where}\;\; k=2,\dots,n.
\end{equation*}

\begin{figure}[h]
\centerline{\includegraphics[width=10cm]       
    {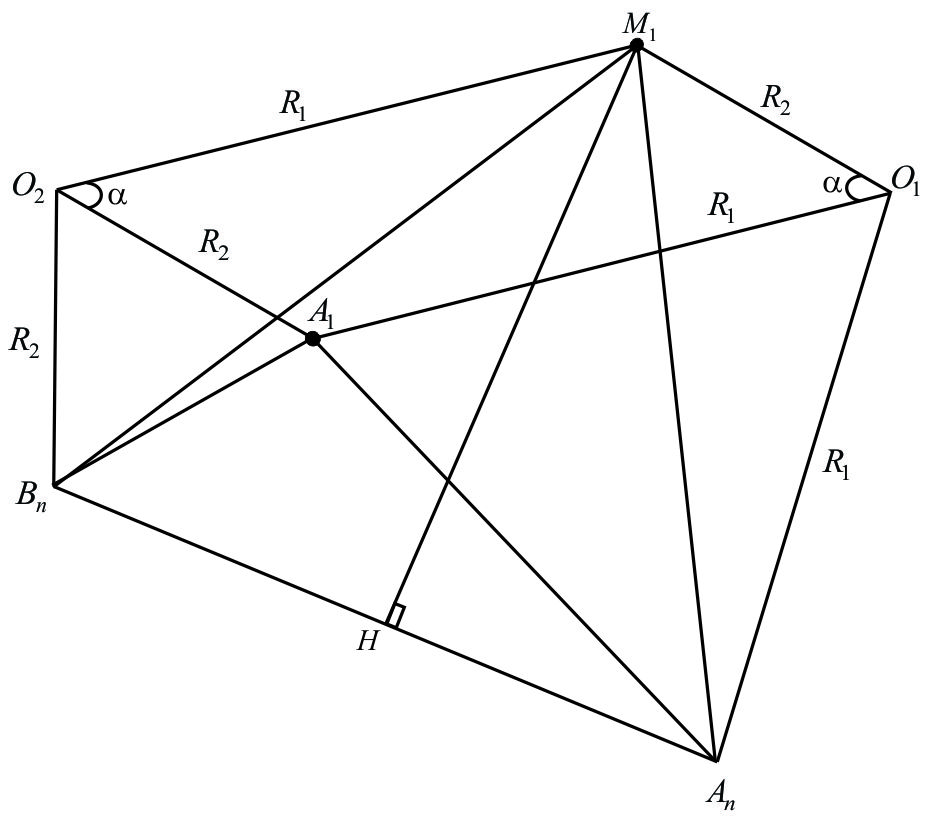} }
\caption{}
\label{fig:9}
\end{figure}

Draw the altitude $M_1H$ of the triangle $A_nM_1B_n$, see Fig.~\ref{fig:9}. The triangle $A_nM_1B_n$ is isosceles $M_1A_n=M_1B_n$. Because
\begin{align*}
    M_1O_1=O_2B_n & =R_2, \\
    M_1O_2=O_1A_n & =R_1;
\end{align*}
so
$$  \vartriangle\!A_nO_1M_1=\vartriangle\!M_1O_2B_n.        $$
Because of $\angle M_1O_1A_1\equiv\alpha$, then
\begin{multline*}
    \angle A_nA_1B_n=2\pi-\big(\angle O_1A_1O_2+\angle O_1A_1A_n+\angle O_2A_1B_n\big) \\
    =\alpha+\frac{2\pi}{n}=\angle A_nO_1M_1=\angle B_nO_2M_1. 
\end{multline*}

The sides $A_1A_n$, $A_1B_n$ are proportional to $R_1$, $R_2$, respectively. Therefore,
$$  \vartriangle\!A_nO_1M_1\sim \vartriangle\!A_nA_1B_n,        $$
so
$$  \angle B_nA_nA_1=\angle O_1A_nM_1.      $$
Then,
\begin{multline*}
    \angle A_nM_1B_n=2\cdot\angle A_nM_1H \\
    =2\Big(\frac{\pi}{2}-\angle B_nA_nM_1\Big)=\pi-2\big(\angle B_nA_nA_1+\angle M_1A_nA_1\big) \\
    =\pi-2\big(\angle O_1A_nM_1+\angle M_1A_nA_1\big)=\pi-2\angle O_1A_nA_1=\angle A_1O_1A_n=\frac{2\pi}{n}\,.
\end{multline*}
For the position of the point $M_1$, we have:
$$  HM_1=\frac{1}{2}\,A_nB_n\cot\frac{\pi}{n}\,.        $$
So the position of the $M_1$ depends only on the line segment $A_nB_n$ and is independent of the location of $A_1$. Moreover, the obtained expression for $HM_1$ also shows that $M_1$ is the center of the regular $n$-gon formed on the side $A_nB_n$.
\ \hfill $\square$
\vskip+0.2cm

For the triangle $\vartriangle\!A_2M_1B_2$ (see Fig.~\ref{fig:7} and Fig.~\ref{fig:8}), analogically we get:
$$  \angle A_2M_1B_2=\frac{2\pi}{n}\,.       $$
i.e.
\begin{equation*}
    \angle A_nM_1B_n=\angle A_2M_1B_2.
\end{equation*}
In the same manner, if we consider the triangles $A_{n-1}M_1B_{n-1}$ and $A_3M_1B_3$, we have
$$  \angle A_{n-1}M_1B_{n-1}=\angle A_3M_1B_3=\angle A_1O_1A_{n-1},     $$
so:

\begin{theorem}\label{th:7}
In two regular $n$-gons $A_1A_2\cdots A_n$ and $A_1B_2\cdots B_n$ with the shared vertex $A_1$ and opposite orientation, take the points $D_1$, $D_2$ on the circumcircles of the polygons, which are diametrically opposed of the shared vertex. Then, for the midpoint $M_1$ of the line segment $D_1D_2$ holds:
$$  \angle A_kM_1B_k=\angle A_{n+2-k}M_1B_{n+2-k}=\frac{2\pi}{n}\,(k-1), \;\;\text{where}\;\; k=2,\dots,n.    $$
\end{theorem}

\bigskip
\bigskip

\section{Conclusion}
\label{sec:9}

\bigskip

In the present paper general method is given -- how to fix two regular polygons in the plane at equal distances from the corresponding vertices of each polygons. The result is trivial for the congruent $n$-gons, but for two non-congruent regular $n$-gons we get an unexpected result, which gives the new theorems in Euclidian geometry. For the existence of the point $M$ two conditions must be satisfied:
\begin{enumerate}
\item[{\rm I.}] $M=\Omega_1(O_2,R_1)\cap\Omega_2(O_1,R_2)$, where $\Omega_1$, $\Omega_2$ -- are circumscibed circles, $R_1$, $R_2$ -- circumradii and $O_1$, $O_2$ -- centroids of the corresponding regular polygons.

\item[{\rm II.}] One pair of the distances from the $M$ must be equal $(t_1=d_1)$, which is obtained by rotation of one of the polygons.
\end{enumerate}

In case of two $n$-gons with the shared vertex, the second condition is satisfied automatically, so we get two points $M_1$ and $M_2$. The properties of these points are investigated. We have established that $M_1$ is the point which satisfies the Bottema theorem $(n=4)$ condition, from where we obtain -- generalized Bottema theorem for any two regular polygons.

\bigskip

\bigskip

\bigskip

DEPARTMENT OF MATHEMATICS

GEORGIAN-AMERICAN HIGH SCHOOL

18 CHKONDIDELI STR., TBILISI 0180, GEORGIA

\textit{E-mail address:} \texttt{mathmamuka@gmail.com}


\bigskip
\bigskip

\begin{thebibliography}{99}

\bigskip

\bibitem{1} M. Meskhishvili, Cyclic averages of regular polygons and platonic solids, Communications in Mathematics and Applications \textbf{11} (2020), no. 3, 335 -- 355, doi: 10.26713/cma.v11i3.1420.


\bibitem{2} M. Meskhishvili, Cyclic averages of regular polygonal distances, \emph{International Journal of Geometry}\textbf{10} (2021), no. 1, 58 -- 65.

\bibitem{3} N. N. Giang, A new proof and some generalizations of the Bottema theorem, \emph{International Journal of Computer Discovered Mathematics (IJCDM)} \textbf{3} (2018),  49 -- 54.

\bibitem{4} R. Seroul, \textit{Newton--Girard Formulas}, In Programming for Mathematicians, pp. 278--279, Springer-Verlag, 2000.

\bibitem{5} A. Shriki, Back to treasure island, \emph{The Mathematics Teacher} \textbf{104} (2011), no. 9, 658 -- 664, doi: 10.5951/MT.104.9.0658.

\bibitem{6} J. P. Tignol, \textit{Galois' Theory of Algebraic Equations}, World Sciences Publications, NJ, 2004.


\end{thebibliography}
\end{document}